\title{Mixed Precision Iterative Refinement with Adaptive Precision Sparse Approximate Inverse Preconditioning\thanks{Version of 
{\today}. Work supported by ERC Starting Grant No. 101075632, SVV grant SVV-2023-260711 and the Exascale Computing Project (17-SC-20-SC), a collaborative
effort of the U.S. Department of Energy Office of Science and the National Nuclear Security Administration. 
}}

\author{Noaman Khan\footnotemark[2]
        \and Erin Carson\footnotemark[2]}
\date{}
\documentclass{article}
\usepackage[utf8]{inputenc}
\usepackage{times}
\usepackage{amsmath}
\usepackage{mathtools}
\usepackage{amssymb}
\usepackage{graphicx}
\usepackage{epstopdf}
\usepackage{color}
\usepackage{graphicx}
\usepackage{lipsum}
\usepackage{graphicx}
\usepackage{subcaption}
\usepackage{hyperref}
\usepackage{algpseudocode}
\usepackage[T1]{fontenc}
\usepackage{hyphenat}
\usepackage[section]{placeins}
\usepackage{multirow,array,graphics}
\usepackage{soul}

\newcommand{\cond}{\text{cond}}

\newcommand{\Pc}{P}

\usepackage[section]{algorithm}

\begin{document}
\maketitle

\begin{abstract}
Hardware trends have motivated the development of mixed precision algorithms in numerical linear algebra, which aim to decrease runtime while maintaining acceptable accuracy.  One recent development is the development of an adaptive precision sparse matrix-vector produce routine, which may be used
to accelerate the solution of sparse linear systems by iterative methods. This approach is also applicable to the application of inexact preconditioners, such as sparse approximate inverse preconditioners 
used in Krylov subspace methods. 

In this work, we develop an adaptive precision sparse approximate inverse preconditioner and demonstrate its use within a five-precision GMRES-based iterative refinement method. We call this algorithm variant BSPAI-GMRES-IR. We then analyze the conditions for the convergence of BSPAI-GMRES-IR, and determine settings under which BSPAI-GMRES-IR will produce similar backward and forward errors as the existing SPAI-GMRES-IR method, the latter of which does not use adaptive precision in preconditioning. Our numerical experiments show that this approach can potentially lead to a reduction in the cost of storing and applying sparse approximate inverse preconditioners, although a significant reduction in cost may comes at the expense of increasing the number of GMRES iterations required for convergence. 

\end{abstract}

\section{Introduction}
We consider the problem of solving large, sparse linear systems $Ax=b$ using iterative methods, where $A$ is a nonsingular $n\times n$ matrix. In recent years, the emergence of low precision, such as half precision, on modern hardware has received renewed attention. 
Lower precision has many benefits, including a reduction in computation, storage, and data movement costs. However, with fewer bits, we have greater round off error and a smaller range of re-presentable numbers. This has motivated the development of mixed precision algorithms, in which lower and higher precisions are used selectively in order to improve the performance, memory, and energy consumption without sacrificing accuracy; for details, see the recent surveys \cite{abdelfattah2021survey, higham2022mixed}. 

Iterative refinement (IR) is a long-standing technique for iteratively improving the solution to a linear system. The idea of iterative refinement is to first compute an initial solution $x_0$ to $Ax=b$ , often using a direct solver like $LU$ factorization. The refinement steps in iterative refinement consist of computing the residual $r_i=b-Ax_i$, solving the correction equation $Ad_i=r_i$, and updating the solution ${x}_{i+1}=x_i+d_i$. In the case that $LU$ factorization is used for computing the initial solution, the LU factors can be reused for solving the correction term $d_i$. This is what we refer to as ``standard IR'' (SIR).

Iterative refinement was originally proposed by Wilkinson in  1948, who suggested performing all the computations in a working precision denoted by $u$ except the residual computation in precision $u^2$. This variant has been analyzed by Wilkinson \cite{wilk63} and Moler \cite{mole67}. In 1977, Jankowski and Wo\'{z}niakowski \cite{jawo77} and Skeel \cite{skee80} introduced fixed precision iterative refinement, performing all computations in precision $u$.
Langou et al. in 2006 used single precision in the computation of the $LU$ factorization, which can be as twice faster as double precision, and a working precision in other parts of the computation \cite{lllk06}. 

The availability of half precision in modern GPUs motivated the development of iterative refinement which uses three or more hardware precisions. Carson and Higham in 2018 proposed an iterative refinement scheme that uses three different precisions $u_f$, $u$, and $u_r$, which denote the factorization, working, and the residual precisions respectively; for an explanation see \cite{carson2018accelerating}. The authors also proposed a fourth precision, called the ``effective precision'', denoted by $u_s$, which allows for general solvers to be used for the correction term $d_i$. 
For example, in standard iterative refinement, the LU factors computed in precision $u_f$ results in $u_s = u_f$. 
With $u_f \geq u$ and $u_r \leq u^2$, then the relative forward and backward errors will converge to level $u$ when $\kappa_\infty(A)\leq u_f^{-1}$, where $\kappa_\infty(A)=\Vert A^{-1}\Vert_\infty \Vert A\Vert_\infty$ denotes the infinity-norm condition number of $A$. 

In \cite{carson2017new}, the authors develop a GMRES-based iterative refinement algorithm (GMRES-IR) which uses the computed $LU$ factors as preconditioners within GMRES to solve for the correction in each refinement step. Under the assumption that GMRES is executed in the working precision $u$, with matrix vector product and with preconditioned matrix computed in double the working precision, $u_s = u$, and thus GMRES-IR is guaranteed to produce forward and backward errors to the working precision for more ill-conditioned problems than standard iterative refinement. 
Assuming that $u_f \geq u$ and $u_r \leq u^2$ a relative forward and backward errors to the level $u$ is obtained for $\kappa_\infty(A)\leq u^{-1/2}u_f^{-1}$.  

From a performance perspective, the requirement that the preconditioned matrix is applied in double the working precision is not attractive.

 In 2021, Amestoy et al. \cite{amestoy2023combining} proposed and analyzed a five-precision variant of GMRES-IR which, in addition to the working precision $u$, factorization precision $u_f$, and residual precision $u_r$, added two more precisions, namely $u_g$ for the working precision within GMRES and $u_p$ for precision in which the preconditioned matrix is applied to a vector within GMRES. The variant with setting $u=u_g=u_p$ is used commonly in practice,  although it is guaranteed to converge for a smaller range of condition numbers than the algorithm in \cite{carson2018accelerating}. Again assuming $u_f \geq u$ and $u_r \leq u^2$, the relative forward and backward error to the level working precision is obtained for the matrices having $\kappa_\infty(A) \leq u^{-1/3}u_f^{-2/3}$, although this restriction is likely overly pessimistic in practice.
 
Most existing analyses of GMRES-based iterative refinement schemes assume that an $LU$ factorization is computed for use as a left preconditioner within GMRES in each refinement step. But when $A$ is very sparse, the performance of this approach may not be attractive since the $LU$ factorization of $A$ may have considerable fill-in. In practice, inexact preconditioners are often used, such as incomplete $LU$ factorizations or sparse approximate inverses (SPAI). Using SPAI has an advantage because it is, in theory, highly parallelizable, as each column can be computed independently, and its application involves only a  sparse matrix-vector product (SpMV). In \cite{carson2022mixed}, the authors propose a new variant called SPAI-GMRES-IR which, instead of $LU$ factors, uses a sparse approximate inverse preconditioner (computed in a precision $u_f$ with a given accuracy threshold $\varepsilon$, which controls the residual in each column) as a preconditioner within five-precision GMRES-IR. The analysis of SPAI-GMRES-IR shows that as long as $\varepsilon$ and $u_f$ satisfy the constraints 
\[
u_f\cond_2(A^T) \lesssim \varepsilon \lesssim u^{-1/2}\kappa_\infty(A)^{-1/2},
\]
then the constraints on condition number for forward error and backward error to converge are the same as for five-precision GMRES-IR with the full $LU$ factors, although it is clear that convergence of the GMRES solves may be slower.

 In 2022, Graillat et al. proposed an adaptive, mixed precision algorithm for computing sparse matrix-vector products that adaptively selects the precision in which each matrix element is stored and applied by splitting them into buckets based on their magnitude and then using progressively lower precisions for the buckets with smaller elements  \cite{graillat2022adaptive}.  
 
 In this work, we apply the idea proposed in \cite{graillat2022adaptive} to the application of the computed SPAI $M$ within SPAI-GMRES-IR. We call this approach BSPAI-GMRES-IR, where the `B' stands for `bucketed'; the components of $M$ are split into different buckets, with a different precision associated with  each bucket. In Section \ref{sec:background} we give background on SPAI preconditioners and the adaptive precision sparse matrix-vector product approach in \cite{graillat2022adaptive}, and discuss bucketed SPAI and recent related approaches. In Section \ref{sec:BSPAI}, we analyze under which conditions the BSPAI-GMRES-IR will converge and bound the forward and backward errors. In Section \ref{sec:numexp} we perform a set of numerical experiments which illustrate the behavior of BSPAI-GMRES-IR. In Section \ref{sec:conclusion} we conclude and discuss future work.

\section{Background}
\label{sec:background}
\subsection{Notation}
First we mention some notation which will be used in rest of the text. Important for us will be the condition numbers. For a given matrix $A$, and a vector $x$, and a norm $p$, we define 
\begin{equation*}
    \kappa_p(A) = {\Vert A^{-1}\Vert}_{p} {\Vert A\Vert}_{p},\hspace{1mm} 
    {\cond}_{p}(A) = {\Vert |A^{-1}||A|\Vert}_{p},\hspace{1mm} 
    {\cond}_{p}(A,x) = \frac{{\Vert |A^{-1}||A||x|\Vert}_{p}}{{\Vert x \Vert}_{p}},
\end{equation*}
where $|A|=(|{a}_{ij}|)$. In case $p$ is not specified we assume the norm to be infinity. For unit roundoffs we will use the notation $u$ and subscripts on $u$ to distinguish various precisions. For rounding error analysis, we will use the notation 
\begin{equation*}
    \gamma_k = \frac{ku}{1-ku}, \quad \tilde\gamma_k=\frac{cku}{1-cku},
\end{equation*}
where $c$ is a small constant independent of problem dimension. A superscript on $\gamma$ indicates that the corresponding $u$ has that superscript as a subscript; for example, $\gamma_k^f = ku_f/(1-ku_f)$. The quantities computed in finite precision will be denoted by hats. 
 
\subsection{Sparse Approximate Inverse Preconditioners}

Sparse approximate inverse preconditioning is based on the idea of explicitly constructing a matrix $M\approx A^{-1}$. Although SPAI is a general algebraic preconditioning technique and is thus not expected to be effective for every problem, the use of SPAI-type  preconditioners within Krylov subspace methods has the advantage that the application of the preconditioner involves only matrix-vector products, unlike, e.g., LU-based preconditioners which require two triangular solves. 

There are many potential techniques for computing a sparse approximate inverse $M$; see the survey \cite{benzi1999comparative, benzi2002preconditioning}. 
A popular approach based on Frobenius norm minimization produces a  sparse approximate inverse in unfactored form  (i.e., a single matrix $M$), in which $M$ is computed as the solution to $\min_{\mathcal{J}\in\mathcal{S}} \Vert I-AM\Vert_F$, where $\mathcal{J}\in \mathbb{B}^{n\times n}$ is a prescribed binary sparsity pattern in the set of all possible binary sparsity patterns $\mathcal{S}\in \mathbb{B}^{n\times n}$. The benefit is that we can decouple this minimization problem as
\begin{equation}
\min_{\mathcal{J}\in\mathcal{S}} \Vert I-AM\Vert_F^2 = \sum_{k=1}^n \min_{\mathcal{J}_k\in\mathcal{S}_k} \Vert e_k-Am_k\Vert_2^2,
\label{eq:fronormmin}
\end{equation}
where $\mathcal{J}_k$, $m_k$, and $e_k$ represent the $k$th columns of $\mathcal{J}$, $M$, and $I$, respectively.
The computed $M$ is then reduced to solving a linear least squares problem for each column $m_k$ of $M$. From a performance point of view, the benefit is that these linear least squares problems are solved independently and in parallel.

Early works based on this approach used a fixed prescribed sparsity pattern $\mathcal{J}$. The set $\mathcal{J}_k$ extracts column indices of $A$ that are relevant for solving for a column $m_k$. The nonzero rows of the submatrix $A(:, \mathcal{J}_k)$ are represented by the so-called ``shadow'' of $\mathcal{J}_k$, 
\[
\mathcal{I}_k = \left\{ i\in\{1,\ldots, n\}: \sum_{j\in \mathcal{J}_k } |a_{ij}|\neq 0\right\},
\]
where $a_{ij}$ is the $(i,j)$ entry of $A$. Thus each term in the summation on the right in \eqref{eq:fronormmin} can be reduced to 
\begin{equation}
\min_{\mathcal{J}(\bar{m}_k) = \mathcal{J}_k} \Vert \bar{e}_k - \bar{A}_k \bar{m}_k \Vert_2,
\label{redprob}
\end{equation}
where $\bar{A}_k = A(\mathcal{I}_k, \mathcal{J}_k)\in \mathbb{R}^{|\mathcal{I}_k|,|\mathcal{J}_k|}$, $\bar{m}_k = m_k(\mathcal{J}_k)\in \mathbb{R}^{|\mathcal{J}_k|}$, $\bar{e}_k = e_k(\mathcal{I}_k)\in\mathbb{R}^{|\mathcal{I}_k|}$, and $\mathcal{J}(\bar{m}_k)$ is the binary sparsity pattern of $\bar{m}_k$. This results in small least squares problems which can be solved directly, for example, via QR factorization.

The deficiency of this approach is that it is hard to predict a sparsity pattern a priori that will ensure an effective preconditioner. Mostly common choices used are the sparsity pattern of $A$, $A^T$, or a power of a sparsified $A$, although generally its not guaranteed that the preconditioner produced will be effective. For overcoming this difficulty, many authors proposed iterative approaches. In one such approach, one starts with an initial sparsity pattern and adds nonzeros to this pattern until $\Vert e_k - Am_k\Vert_2\leq \varepsilon$ becomes true for some threshold $\varepsilon$ or the maximum number of nonzeros has been reached.  
For a more detailed explanation of this type algorithm, see, e.g.,  the work by Cosgrove et al. \cite{cosgrove1992approximate}, Grote and Huckle \cite{grote1997parallel}, and Gould and Scott \cite{gould1998sparse}.

The most successful among these algorithms is that of Grote and Huckle \cite{grote1997parallel} which is commonly used to compute a SPAI preconditioner \cite{benzi1999comparative}, and which we use in the present work. To overcome the difficulty of choosing the sparsity pattern a priori for a resulting effective preconditioner, the authors in \cite{grote1997parallel} proposed an adaptive approach that dynamically determines the most beneficial nonzero indices to include. 
Algorithm \ref{alg:spai} is one specific variant of Grote and Huckle's algorithm, which is taken from \cite[Algorithm 4]{sedlacek2012sparse}. 

The algorithm requires an input matrix $A$,  
$\mathcal{J}$ as the initial binary sparsity pattern, $\varepsilon$ as the convergence tolerance, $\alpha$, for the maximum number of iterations for each column, and $\beta$, for the maximum number of nonzeros added to the pattern in each iteration.

The algorithm for each column solves the linear least squares problem \eqref{redprob} for a given initial sparsity pattern $\mathcal{J}$  and computes the residual $\bar{s}_k$ (lines \ref{solv1}-\ref{solvn}). This column is considered finished when the 2-norm of the residual is less than the threshold $\varepsilon$. Otherwise, we continue adding entries to $\mathcal{J}$.

We construct an index set $\mathcal{L}_k$ in line \ref{ell} which contain the nonzeros entries in $\bar{s}_k$. 

From the index set $\mathcal{L}_k$, for every element $\ell$ we go through that $\ell$th row of $A$ and choose the column indices of the nonzero entries for which we define a set name $\tilde{\mathcal{J}}_k$ which are not $\mathcal{J}_k$. The set $\tilde{\mathcal{J}}_k$ is the union of the sets $\mathcal{N}_\ell$ which contain the potential indices that can be added to  $\mathcal{J}_k$, out of which we select only a subset of the ``most important'' indices.

There are many ways to determine which indices are most important. Grote and Huckle's technique considers a univariate minimization problem,  through which the quantity $\rho_{jk}$ computed in line \ref{rhojk} gives a measure of the 2-norm of the new residual if index $j$ is added to $\mathcal{J}_k$. A well-known heuristic (see, e.g., \cite{barnard1999mpi}) is to mark indices as ``acceptable'' if their $\rho_{jk}$ is less than the arithmetic mean $\bar{\rho}_k$ over all $j$. Then we choose up to $\beta$ of the best (smallest $\rho_{jk}$) indices acceptable to add (lines \ref{add1}-\ref{addn}) in each of the $\alpha$ iterations.

In line \ref{qrfact} there is no need to recompute the QR factorization fully in each step; the factorization can be updated by using the QR factorization computed in the previous step and the entries added to $\bar{A}_k$; see \cite[Eqns.~(14) and (15)]{grote1997parallel}. Typical values for the parameters are $\varepsilon\in [0.1,0.5]$, $\alpha\in\{1,\ldots,5\}$, and $\beta\in\{3,\ldots,8\}$ \cite[Section 3.1.3]{sedlacek2012sparse}.

In SPAI, although each column can theoretically be computed in parallel, the construction is often costly, specially for large-scale problems; see, e.g., \cite{barnard1999mpi, benzi1999comparative, chow2001parallel, gao2021thread}. SPAI  memory requirements scale quadratically and the computational cost scales cubically in the number of nonzeros per row \cite{gao2021thread}. Thus applying the bucketing idea to sparse approximate inverse preconditioner in which low precision is used for the buckets containing elements of smaller magnitude has the potential to significantly reduce this cost. 
For modern hardware like GPUs, the construction of efficient sparse approximate inverse computations has been the subject of much recent work; see, e.g., \cite{gao2017gpu,lukash2012sparse,dehnavi2012parallel,he2020efficient}.

\begin{algorithm}[h]
\caption{Variant of sparse approximate inverse (SPAI) construction (\cite[Algorithm 4]{sedlacek2012sparse})}\label{alg:spai}
\begin{algorithmic}[1]
\Require {$A \in\mathbb{R}^{n\times n}$, $\mathcal{J} \in \mathbb{B}^{n\times n}$, $\alpha \geq 0$, $\beta \geq 0$, $\varepsilon >0$}
\Ensure {Right preconditioner $M \approx A^{-1}$, $M\in\mathbb{R}^{n\times n}$}
\For{$k=1$ to $n$}
    \State {$e_k=I(:,k)$}
    \State {$\mathcal{J}_k=\mathcal{J}(:,k)$}
    \For{step $=0$ to $\alpha$}
        \State {$\mathcal{I}_k=\bigg\{i \in \left\{1, \ldots ,n\right\}\colon \sum_{j \in \mathcal{J}_k} \left|a_{ij}\right| \neq 0 \bigg\}$} \label{solv1}
        \State {$\bar{A}_{k}=A(\mathcal{I}_k,\mathcal{J}_k)$}
        \State {$\bar{e}_{k}=e_k(\mathcal{I}_k)$}
        \State {Compute QR factorization $\bar{A}_{k} = \bar{Q} \bar{R}$} \label{qrfact}
        \State {$\bar{m}_{k}=\bar{R}^{-1}\bar{Q}^{T}\bar{e}_{k}$}
        \State {$\bar{s}_{k}= \bar{A}_k\bar{m}_{k}-\bar{e}_{k}$} \label{solvn}
        \If{$\Vert \bar{s}_{k}\Vert_{2}$} $\leq$ $\varepsilon$
            \State {break}
        \EndIf
        \State {$\mathcal{L}_k=\mathcal{I}_k$ $\cup$ $\left\{ k \right\}$} \label{ell}
        \For {$\ell$ $\in$ $\mathcal{L}_k$}
            \State {$\mathcal{N}_\ell=\{j \colon {a}_{\ell j} \neq 0 \}$} \label{enn}
        \EndFor
        \State {$\tilde{\mathcal{J}}_k =\bigcup_{\ell \in \mathcal{L}_k} \mathcal{N}_\ell$}
        \State {$\tilde{\rho}_{k}=0$}
        \For {$j$ $\in$ $\tilde{\mathcal{J}}_k$}
            \State {${\rho}_{jk}=\left(\Vert \bar{s}_{k}\Vert^{2}_{2}-\frac{[\bar{r}^{T}_{k}A_j(\mathcal{I}_k)]^2}{\Vert A_j(\mathcal{I}_k)\Vert^{2}_{2}}\right)^\frac{1}{2}$} \label{rhojk}
            \State {$\tilde{\rho}_{k}=\tilde{\rho}_{k}$ + $\tilde{\rho}_{jk}$}
        \EndFor
        \State {$\tilde{\rho}_{k}=\frac{\tilde{\rho}_{k}}{|\tilde{\mathcal{J}}_{k}|}$}
        \For{idx $=1$ to $\beta$}
            \State {$j={argmin}_{j \in \tilde{\mathcal{J}}_{k}}$ $\rho_{jk}$} \label{add1}
            \State {$\mathcal{J}_k={\mathcal{J}}_{k}$ $\cup$ $\left\{j \colon \rho_{jk} \leq \tilde{\rho}_{k}\right\}$}
            \State{$\tilde{\mathcal{J}}_{k}=\tilde{\mathcal{J}}_{k}$\textbackslash $\left\{j\right\}$} \label{addn}
        \EndFor
    \EndFor
    \State {$m_k(\mathcal{J}_k)=\bar{m}_{k}$}
\EndFor
\end{algorithmic}
\end{algorithm}

\subsection{Adaptive Precision Sparse Matrix-Vector Products}
\label{sec:adaptivespmv}
As mentioned, with the emergence of low precision arithmetic,  such as half precision fp16 or bfloat16 on modern computers, mixed precision algorithms in numerical linear algebra have received renewed attention. Many variants of mixed precision algorithms have been recently proposed; see, for example, the works \cite{markidis2018nvidia, ahmad2019data, blanchard2020mixed, mukunoki2020dgemm} on matrix multiplication. The works \cite{carson2017new, carson2018accelerating, aadg16, h:21} proposed mixed precision iterative refinement methods based on preconditioned Krylov subspace methods. The authors in \cite{higham2019new} proposed a general preconditioning technique based on a low-rank approximation of the error. 

A particularly fruitful idea is the concept of \emph{adaptive} precision algorithms, in which the precisions used need not be determined a priori, but are instead dynamically set based on the data involved in the computation and perhaps some user-specified accuracy constraints. Often, the precisions chosen are proportional to importance of the data, which is inherently application dependent. For example, the authors in \cite{anzt2019adaptive, flegar2021adaptive} introduced an adaptive precision block Jacobi preconditioner with idea of choosing the precision of each block based on its condition number. Amestoy et al. \cite{amestoy2021mixed} introduced mixed precision block low rank compression that partitions a low rank matrix into several low-rank components of decreasing norm and stores each of them in a correspondingly decreasing precision. Ahmad et al. \cite{ahmad2019data} introduced an algorithm for sparse matrix-vector products that switches the elements in the range of $[-1, 1]$ to single precision while keeping the other elements in double precision. The authors in \cite{diffenderfer2021qdot} develop a ``quantized'' dot product algorithm, adapting the precision of each vector element based on its exponent. 

In recent work, which is the focus of the present paper, Graillat et al. \cite{graillat2022adaptive} develop an adaptive precision sparse matrix-vector product algorithm with the idea of adapting the precision of each matrix element based on its magnitude. The elements of the matrix are split into different buckets and different precisions are used to store and compute with elements in each bucket. Buckets with smaller elements are stored in lower precision. This approach is used to apply the matrix $A$ to a vector within GMRES-IR with Jacobi preconditioning. We now give an overview of the results of \cite{graillat2022adaptive}.

\begin{algorithm}[ht]
\caption{Bucketing Procedure}\label{alg:bspai}
 \begin{algorithmic}[1]
\Require {$A \in\mathbb{R}^{n\times n}$, precisions $u_1,...,u_q$,  target accuracy $\epsilon$}
\Ensure {Buckets $B_{ik}$, for $i\in\{1,\ldots,n\}$, $k\in \{1,\ldots,q\}$}
\For{$i=1:n$}
  \State {$J_i$ = $\text{find}(A_{i,:})$}
     \For{$i=1:q$}
         \State{$B_{ik}$ = $\{\}$}
    \EndFor
    \For{$l=1$ to $|J_i|$}
         \State{$j=J_i(l);$}
         \If{$(|a_{ij}| > \epsilon \|A\|/u_2$}
             \State{$B_{i1}=\{ B_{i1}, j\}$}
         \ElsIf{$(|a_{ij}| \leq \epsilon \|A\|/u_q$}
            \State{$B_{iq}=\{B_{iq}, j\}$}
         \Else   
            \For{$k=2:q-1$}
                \If{$(|a_{ij}| > \epsilon \|A\|/u_{k+1})$ and $(|a_{ij}|\leq \epsilon \|A\|/u_k$}
                    \State{$B_{ik}=\{B_{ik}, j\}$}
                \EndIf
            \EndFor
        \EndIf
    \EndFor         
\EndFor
\end{algorithmic}
\end{algorithm}

For matrix-vector products in a uniform precision, the Oettli-Prager\cite[Thm.7.3]{high:ASNA2}, \cite{oettli1964compatibility} and Rigal-Gaches\cite[Thm.7.1]{high:ASNA2}, \cite{rigal1967compatibility} theorems give the formula for normwise backward error,  

\begin{equation}
\varepsilon_{nw}=\min\{\varepsilon:\hat{y} = (A+\Delta A)x, \, \| \Delta A \|\leq \varepsilon \|A\| \ \} = \frac{\|\hat{y}-y\|}{\|y\| \|x\|}.
\label{nw}
\end{equation}

A bound on the normwise backward error for the uniform precision case is  
\begin{equation*}
\varepsilon_{nw} \leq pu,
\end{equation*}
where $p$ is the maximum number of nonzero elements per row of $A$; see, e.g., \cite[Eqn (2.2)]{graillat2022adaptive}.

The idea of the adaptive precision sparse-matrix vector product  approach of Graillat et al. \cite{graillat2022adaptive} is, for a given set of $q$ precisions i.e $u_1 < u_2 \ldots < u_q$, to split the elements of the matrix $A$ into $q$ buckets based on the magnitude of the elements. Using this approach splits the nonzeros elements in each row $i$ of the computed $M$ into up to $q$ buckets and then computes the partial inner products associated with each bucket in up to $q$ different precisions. The partial inner products are then all summed in precision $u_1$. 

We briefly recall the notation, algorithm, and key points of the error analysis given in \cite{graillat2022adaptive}.
Let $J_i$ denote the set of column indices of the nonzero elements in row $i$ of $A$. Each row $i$ of the matrix $A$ will be partitioned into the $q$ buckets $B_{ik}\subset [1,n]$ for $k=1:q$. How we define the buckets will affect the resulting normwise (or componentwise) backward error. Assume that we want to construct the buckets ${B}_{ik}$ in such a way that the backward error obtained is at most of order $O(\epsilon)$, where $\epsilon$ is the user defined target accuracy with $\epsilon$ $\geq$ $u_1$. We can define the buckets as 
\[
B_{ik} = \{ j\in J_i : |a_{ij}| \in P_{ik}\},
\]
with 
\begin{equation}
{P}_{ik}= 
\begin{cases}
\bigl(\epsilon\lVert A \rVert/u_2, \, \, +\infty \bigl) & \text{for}\, k=1 \\
\bigl(\epsilon\lVert A \rVert /{u}_{k+1},\,\, \epsilon\lVert A \rVert/{u_k}\bigr] & \text{for} \, k=2:q-1 \\
\left[0, \, \epsilon\lVert A \rVert/u_q\right] & \text{for} \, k=q 
\label{nint}
\end{cases}.
\end{equation}
The procedure for placing elements of a matrix $A$ into buckets according to this rule is given in Algorithm \ref{alg:bspai}.

The partial inner product $y_i^{(k)} = \sum_{j\in B_{ik}} a_{ij}x_j$ associated with bucket $B_{ik}$ is computed in precision $u_k$, and all partial inner products are accumulated in precision $u_1$ (the highest precision). This procedure is given in Algorithm \ref{alg:BSPMV}. Theorem 3.1 in \cite{graillat2022adaptive} states that if $y=Ax$ is computed using this approach, then we have 
\begin{equation}
    \varepsilon_{nw} \leq (q-1) u_1 + c\epsilon,
\end{equation}
where 
\begin{equation}
 c= (1+ (q-1)u_1 )+ {\max}_{i}{\sum}_{k=1}^{q}{p}_{ik}^{2}{(1+u_k)}^{2},
 \label{ccw}
 \end{equation}
and $p_{ik}$ is the number of elements in $B_{ik}$.

We note that Graillat et al. also provide different bucketing strategies that give guaranteed bounds on the componentwise backward error. The drawback of these is that the bucketing scheme depends on the values in the vector $x$ to be multiplied, and thus the bucketing would need to be redone for each matrix-vector product encountered. Thus for practical reasons we restrict ourselves to the variant which provides normwise error bounds. 
    
\begin{algorithm}[ht]
\caption{Adaptive Precision Matrix-Vector Product \cite[Alg. 3.1]{graillat2022adaptive}}\label{alg:BSPMV}
 \begin{algorithmic}[1]
\Require {$A \in\mathbb{R}^{m\times n}$, $x$, ${B}_{ik}$}
\Ensure {$y=Ax$}
\For{$i=1:m$}
     \For{$k=1:q$}
         \State{${y}_{i}^{(k)}=0$}
         \For{$j \in {B}_{ik}$}
         \State{${y}_{i}^{(k)}={y}_{i}^{(k)} + {a}_{ij}x_{j}$ in precision $u_k$} 
         \EndFor
    \EndFor     
    \State{$y_{i} = \sum_{k=1}^{k=q} {y}_{i}^{(k)}$ in precision $u_1$}
\EndFor
\end{algorithmic}
\end{algorithm}

\section{GMRES based Iterative Refinement with BSPAI}
\label{sec:BSPAI}

Our approach will be to apply the adaptive precision sparse-matrix vector product described in Section \ref{sec:adaptivespmv} to the application of a sparse approximate inverse preconditioner $M$ computed using \ref{alg:spai} within GMRES-based iterative refinement. The resulting algorithm, which we refer to as BSPAI-GMRES-IR, is given as Algorithm \ref{alg:adapspaigmresir}.

\begin{algorithm}[h]
\caption{GMRES-based Iterative Refinement with Adaptive-Precision SPAI Preconditioning (BSPAI-GMRES-IR)}\label{alg:adapspaigmresir}
\begin{algorithmic}[1]
\Require {$n\times n$ nonsingular matrix $A$ and length-$n$ right-hand side vector $b$, maximum number of refinement steps $i_{max}$, GMRES convergenge tolerance $\tau$, SPAI parameter ${\varepsilon}$, parameters for SPAI into buckets (${\epsilon}$,$u_1,\ldots, u_q$).}
\Ensure {Approximate solution $x_{i+1}$ to $Ax=b$.}
\State{Run Algorithm \ref{alg:spai} on $A^T$ with parameter ${\varepsilon}$ in precision $u_f$ to obtain $M$; set $M = M^T$.}
\State {Compute {$x_0=M b$} in precision $u_f$; store $x_0$ in precision $u$.}
\For {$i = 0:i_{max}-1$}
\State {Compute $r_i=b-Ax_i$ in precision $u_r$; store in precision $u$.}
\State{Solve $M Ad_i = M r_i$ via left-preconditioned GMRES with tolerance $\tau$ in working precision $u_g$, with matrix-vector products with $M$ applied using Alg. \ref{alg:BSPMV} and $A$ computed in precision $u_p$; store $d_i$ in precision $u$.}
\State {Update $x_{i+1}=x_i+d_i$ in precision $u$.}
\EndFor
\end{algorithmic}
\end{algorithm}

Our aim is to derive the conditions under which BSPAI-GMRES-IR (Algorithm \ref{alg:adapspaigmresir}) will converge. 

We can determine the resulting backward and forward errors in GMRES when we use the adaptive precision SpMV to apply the preconditioner $M$ within each GMRES iteration. We will assume here that matrix-vector products with $A$ are computed in precision $u_p$ within GMRES (where we will generally take $u_p=u_g=u$, using the notation of \cite{h:21}). Note that we could, in principle, also use the adaptive precision SpMV to apply $A$ to a vector; extending the analysis to this case is simple and the results will not be significantly different as long as $u_p \approx \epsilon_{bspai}$. We give backward and forward error bounds for GMRES for this case as well below. 

Following \cite{h:21} and \cite{carson2017new}, let $z_j = MA \hat{v}_j$ be computed in each iteration of MGS-GMRES as described above, where $A$ is applied in precision $u_p$ and $M$ is applied using the adaptive precision SpMV approach (Algorithm \ref{alg:BSPMV}). Then 
\begin{align*}
    (A+\Delta A) \hat{v}_j &= \hat{w}_j, \qquad \| \Delta A \|_F\leq \gamma_q^p \|A\|_F \\
    (M + \Delta M) \hat{w}_j &= \hat{z}_j, \qquad \|\Delta M\|_F \leq \left((q-1)u_1 + c\epsilon \right) \|M\|_F.
\end{align*}

Then 
\begin{align*}
    \hat{z}_j &= (M+\Delta M)(A+\Delta A) \hat{v}_j \\
    &\approx (MA + M\Delta A + \Delta MA)\hat{v}_j \\
    &= MA\hat{v}_j + f_j,
\end{align*}
where $f_j = (M\Delta A + \Delta A M)\hat{v}_j$. We can bound the norm of this quantity by 
\begin{align*}
    \| f_j \|_2 &\lesssim \gamma_q^p \|M\|_F \|A\|_F + \left( (q-1) u_1 + c \epsilon\right) \|M\|\|A\|\\
    &\leq (q u_p + (q-1)u_1 + c\epsilon) \|M \|_F \|A\|_F \|\hat{v}_j \|_2.
\end{align*}

This means that we can apply \cite[Theorem 3.1]{h:21} with 
\[
\epsilon_p = \left(qu_p + (q-1)u_1 + c \epsilon \right) \frac{\|M\|_F\|A\|_F}{\|MA\|_F}.
\]

Note also that we must apply the preconditioner $M$ to the right-hand side $\hat{r}_i$. Denoting $s_i = M\hat{r}_i$, the computed $\hat{s}_i$ satisfies
\[
\hat{s}_i = (M+\Delta M) \hat{r}_i = s_i + \Delta M \hat{r}_i. 
\]
We then have 
\begin{align}
\| \hat{s}_i - s_i \|_\infty &\leq \left((q-1)u_1+c\epsilon\right) \|M\|_\infty \|\hat{r}_i \|_\infty \nonumber \\
&\leq \left((q-1)u_1+c\epsilon\right) \kappa_\infty(M) \|s_i \|_\infty. \label{eq:SHMS}
\end{align}

Letting $\tilde{A}=MA$, and assuming we are solving the $n\times n$ linear system $\tilde{A}d_i=\hat{s}_i$, the conclusions of \cite[Theorem 3.1]{h:21} say that for MGS-GMRES in working precision $u_g$, except for products with $\tilde{A}$ which satisfy 
\[
fl(\tilde{A}v) = \tilde{A}v + f, \qquad \|f\|_2 \lesssim \epsilon_p \|\tilde{A}\|_F \|v\|_2,
\]
as long as 
\[
\sigma_\text{min}(\tilde{A}) \gtrsim \Big(k^{1/2}\left(qu_p + (q-1)u_1 + c \epsilon \right) \frac{\|M\|_F\|A\|_F}{\|\tilde{A}\|_F} + \tilde{\gamma}_{kn}^g\Big) \|\tilde{A}\|_F,
\]
then for some step $k\leq n$, the algorithm produces an approximate solution $\hat{d}_i$ satisfying
\begin{align}
    (\tilde{A} + \Delta \tilde{A}) \hat{d}_i &= \hat{s}_i + \Delta \hat{s}_i, \label{eq:BE}\\
    \|\tilde{A}\|_F &\lesssim \Big(k^{1/2}\left(qu_p + (q-1)u_1 + c \epsilon \right) \frac{\|M\|_F\|A\|_F}{\|\tilde{A}\|_F} + \tilde{\gamma}_{kn}^g\Big) \|\tilde{A}\|_F, \label{eq:DTA} \\
    \|\Delta \hat{s}_i\|_2 &\lesssim \tilde{\gamma}_{kn}^g \|\hat{s}_i\|_2 \lesssim n^{1/2}\tilde{\gamma}_{kn}^g \|s_i\|_{\infty}. \label{eq:DSH}
\end{align}

From \eqref{eq:BE}, we can write
\[
s_i - \tilde{A}\hat{d}_i = \Delta \tilde{A} \hat{d}_i - (\hat{s}_i - s_i ) - \Delta \hat{s}_i,
\]
which we can bound using \eqref{eq:DTA}, \eqref{eq:SHMS}, and \eqref{eq:DSH}, giving 
\begin{align*}
    \Vert s_i - \tilde{A}\hat{d}_i \Vert_\infty &\leq \| \Delta \tilde{A}\|_\infty \|\hat{d}_i\|_\infty + \|\hat{s}_i - s_i\|_\infty - \Delta \hat{s}_i \\
    &\leq n \Big(k^{1/2}\left(qu_p + (q-1)u_1 + c \epsilon \right) \frac{\|M\|_F\|A\|_F}{\|\tilde{A}\|_F} + \tilde{\gamma}_{kn}^g\Big) \|\tilde{A}\|_\infty \|\hat{d}_i\|_\infty \\
    &\phantom{\leq}+ \left((q-1)u_1+c\epsilon\right) \kappa_\infty(M) \|s_i \|_\infty + n^{1/2}\tilde{\gamma}_{kn}^g \|s_i\|_{\infty} \\
    &\leq  n \Big(k^{1/2}n\left(qu_p + (q-1)u_1 + c \epsilon \right) \kappa_\infty(M) + \tilde{\gamma}_{kn}^g\Big) \|\tilde{A}\|_\infty \|\hat{d}_i\|_\infty \\
    &\phantom{\leq}+ \left((q-1)u_1+c\epsilon\right) \kappa_\infty(M) \|s_i \|_\infty + n^{1/2}\tilde{\gamma}_{kn}^g \|s_i\|_{\infty} \\
    &\leq kn^2 \Big( u_g + \big(qu_p + (q-1)u_1+c\epsilon \big)\kappa_\infty(M) \Big) \Big( \|\tilde{A}\|_\infty \|\hat{d}_i\|_\infty + \|s_i \|_\infty\Big).
\end{align*}

Thus the normwise relative backward error of the system $\tilde{A}\hat{d}_i = s_i$ is bounded by 
\begin{equation}
    \frac{\|s_i - \tilde{A}\hat{d}_i\|_\infty}{\|\tilde{A}\|_\infty \|\hat{d}_i\|_\infty + \|s_i \|_\infty} \lesssim f(n,k)\Big( u_g + \big(qu_p + (q-1)u_1+c\epsilon \big)\kappa_\infty(M) \Big),
    \label{eq:nbe}
\end{equation}
and thus the relative error of the computed $\hat{d}_i$ is bounded by 
\begin{equation}
    \frac{\|d_i -\hat{d}_i\|_\infty}{\|d_i\|_\infty} \lesssim f(n,k)\Big( u_g + \big(qu_p + (q-1)u_1+c\epsilon \big)\kappa_\infty(M) \Big) \kappa_\infty(\tilde{A}),
    \label{eq:fe}
\end{equation}
where $f(n,k) = kn^2$. From \eqref{eq:nbe} and \eqref{eq:fe}, we can say that if $u_1\approx \epsilon \approx u_p$, then the backward and forward errors in MGS-GMRES with adaptive precision SpMV used to apply $M$ will be approximately the same as the case of uniform precision SpMV; see \cite[Eqns. (3.17) - (3.18)]{h:21}. 

We note that in the case where we use the adaptive precision SpMV also in applying the matrix $A$ to a vector within GMRES, the bound for the normwise relative backward error in \eqref{eq:nbe} becomes
\begin{equation}
    \frac{\|s_i - \tilde{A}\hat{d}_i\|_\infty}{\|\tilde{A}\|_\infty \|\hat{d}_i\|_\infty + \|s_i \|_\infty} \lesssim f(n,k)\Big( u_g + \big(2(q-1)u_1+(c_A+c_M)\epsilon \big)\kappa_\infty(M) \Big),
    \label{eq:nbe2}
\end{equation}
where we assume that the same buckets are used for both $M$ and $A$, and $c_A$ and $c_M$ are the values of $c$ in \eqref{ccw} associated with $A$ and $M$, respectively. Similarly, the relative forward error becomes 
\begin{equation}
    \frac{\|d_i -\hat{d}_i\|_\infty}{\|d_i\|_\infty} \lesssim f(n,k)\Big( u_g + \big(2(q-1)u_1+(c_A+c_M)\epsilon \big)\kappa_\infty(M) \Big) \kappa_\infty(\tilde{A}).
    \label{eq:fe2}
\end{equation}
Thus if $u_1\approx \epsilon$, MGS-GMRES with the adaptive precision SpMV used for applying both $M$ and $A$ will produce backward and forward errors similar to the MGS-GMRES variant in \cite{h:21} with the setting $u_p = u_1$.

\section{Numerical Experiments}
\label{sec:numexp}

We perform numerical experiments to evaluate the performance of BSPAI-GMRES-IR by comparing it with SPAI-GMRES-IR in \cite{carson2022mixed}. \emph{We stress that we only expect BSPAI-GMRES-IR to have a clear potential advantage over SPAI-GMRES-IR for the case $u_f=u$.} Otherwise, for example, if $u_f=$ half and $u=$ single, SPAI-GMRES-IR stores the preconditioner entirely in precision $u_f$ but applies it in precision $u$. BSPAI-GMRES-IR, on the other hand, stores the preconditioner in multiple precisions, where we must have $u_1=$ single in order to enable reading effective application precision $\epsilon \approx u$. We also note that this motivates future work in the direction of decoupling the storage and application precisions in adaptive precision sparse matrix-vector products.

All the experiments are performed in \text{MATLAB R2021a}.
The matrices we tested are taken from the SuiteSparse Matrix Collection \cite{davis2011university}. We run the experiments using four precisions which are half, single, double, and quadruple. For properties of these precision, see Table \ref{tab:ieeeprecisions}. 

For half precision, we use the \texttt{chop} library\footnote{\texttt{https://github.com/higham/chop}}. We use MATLAB built-in datatypes for single and double precision and the Advanpix Multiprecision Computing Toolbox for quadruple precision; see \cite{adva-mct}. The code for reproducing the experiments in this paper is available online\footnote{\texttt{https://github.com/Noaman67khan/ADAP-SPAI-GMRES-IR}}.

Matrices used in the experiments along with their key properties are listed in Table \ref{tab:mats}. We set the right-hand side to the vector with equal components and unit 2-norm in all tests. 
For the GMRES tolerance, we set $\tau = 10^{-4}$ in the case working precision is single and $\tau = 10^{-8}$ for the case working precision is double, which responds to roughly the square root of the working precision. These values are set by default used in the previous works by \cite{h:21}, \cite{carson2018accelerating}, and are also used in practical applications. In all invocations, we use the GMRES setting $u_g=u_p=u$, which is commonly used in practice. 

We tested the matrices in Table \ref{tab:mats} with a subset of the settings $(u_f, u, u_r)=$ (double, double, quad), $(u_f, u, u_r)=$  (single, double, quad), and $(u_f, u, u_r)=$ (half, single, double), depending on whether SPAI-GMRES-IR converges with the given precisions and value of $\tau$.
We choose the identity matrix as the initial sparsity pattern for SPAI in all tests. When $A$ has zero entry on the diagonal, this results in a zero column in the SPAI preconditioner, as mentioned in Sedlacek \cite[Section 3.1.2]{sedlacek2012sparse}. Therefore we only choose problems with nonzero entries on the diagonal, but note that this could be remedied by either permuting $A$ or using the initial sparsity pattern of $A$, which, when SPAI is run on $A^T$, guarantees that we obtain a $\Pc$ with nonzero rows \cite[Theorem 3.1]{sedlacek2012sparse}.

In all tests, the matrices are preprocessed with column scaling such that the absolute value of the largest value in every column of $A^T$ is $1$. This one-sided scaling was proposed in \cite[Section 3.2]{carson2020three} to avoid overflow in the computation of $QR$ due to low precision. To be specific, for obtaining $M$, we run SPAI on the scaled $A^T D$ and then set $M = M^T D$, where is the $D$ diagonal scaling matrix. 

In both BSPAI-GMRES-IR and SPAI-GMRES-IR we use the variant in which $u_g=u_p=u$, which is commonly used in practice.  For all tests, we use $\beta=8$, which is in the range suggested by Sedlacek  \cite{sedlacek2012sparse}.
 
For BSPAI-GMRES-IR, when $u$ is double, we use the precisions with $u_1=$ double, $u_2$ = single, $u_3$ = half, and $u_4=1$. When $u$ is single, we use the precisions $u_1$ = single, $u_2$ = half, and $u_3=1$. Note that the choice $u_1=1$ enables the dropping of elements in $M$, as described in \cite[Remark 3.2]{graillat2022adaptive}.

For each linear system and given combination of precisions, we run BSPAI-GMRES-IR with various values of ${\epsilon}\geq u_1$,  and use the same value of $\varepsilon$ for both BSPAI-GMRES-IR and SPAI-GMRES-IR. We report our results in a series of tables. The first column of the table lists the matrix name, and the second column indicates whether we use BSPAI or SPAI and the corresponding parameters. The third column gives the infinity-norm condition number of the preconditioned coefficient matrix. The fourth column gives information about the number of nonzeros and their storage precisions. The first number gives the total number of nonzeros, and the tuple that follows gives information about the precisions: element $i$ in the tuple gives the number of nonzeros stored in precision $u_i$. The fifth column gives the storage cost of the BSPAI preconditioner with mixed precision storage as a percentage of the cost of the SPAI preconditioner with uniform precision storage (the lower the better). The final column gives information about convergence of the iterative refinement process. The first number gives the total number of GMRES iterations over all refinement steps, and element $i$ of the tuple that follows gives the number of GMRES iterations in refinement step $i$. Thus the number of elements in the tuple gives the number of iterative refinement steps required until convergence of the forward and backward errors to the level of the working precision. 

For each setup we form one table with five columns in which first one represent matrices names, second for the preconditioner( SPAI and BSPAI), third for the condition number of the preconditioned system, fourth for the total number of nonzeros (number of nonzeros in each bucket) and the last column is for the information about the number of GMRES-IR refinement steps and GMRES iterations per refinement step.
\begin{table}[t]
\centering
\caption{List of matrices and their properties we used in experiments. All matrices are taken from the SuiteSparse collection \cite{davis2011university}. }
\label{tab:mats}
\begin{tabular}{|c|c|c|c|c|}
\hline
Name & $n$ & $nnz$ & $\kappa_\infty(A)$ & $\text{cond}_2(A^T)$ \\ \hline
\texttt{pores\_3} & 532 & 3474 & $1.2\cdot 10^{6}$ & $1.7\cdot 10^{5}$ \\ \hline
\texttt{bfw782a} & 782 & 7514 & $6.8\cdot 10^{3}$ & $1.3\cdot 10^{3}$ \\ \hline
\texttt{cage5} & 37 & 233 & $2.9\cdot 10^{1}$ & $7.5\cdot 10^{0}$ \\ \hline
\texttt{gre\_115} & 115 & 421 & $1.4\cdot 10^{2}$ & $4.5\cdot 10^{1}$ \\ \hline
\texttt{hor\_131} & 434 & 4182 & $1.5\cdot 10^{5}$ & $3.3\cdot 10^{3}$ \\ \hline
\texttt{sherman4} & 1104 & 3786 & $3.1\cdot 10^3$ & $2.1\cdot 10^3$ \\ \hline
\texttt{steam1} & 240 & 2248 & $3.1\cdot 10^7$ & $2.8\cdot 10^3$ \\ \hline
\texttt{steam3} & 80 & 314 & $7.6\cdot 10^{10}$ & $5.6\cdot 10^3$ \\ \hline
\end{tabular}
\end{table}

\begin{table}[H]
\centering
\caption{Parameters for IEEE floating point precisions. The range denotes the order of magnitude of the largest and smallest positive normalized floating point numbers. }
\label{tab:ieeeprecisions}
\begin{tabular}{l|l|l|l}
Type & Size & Range & Unit Roundoff $u$ \\ \hline
half & 16 bits & $10^{\pm 5}$ & $2^{-11} \approx 4.9\cdot 10^{-4}$ \\
single & 32 bits & $10^{\pm 38}$ & $2^{-24} \approx 6.0\cdot 10^{-8}$ \\
double & 64 bits & $10^{\pm 308}$ & $2^{-53} \approx 1.1\cdot 10^{-16}$ \\
quad & 128 bits & $10^{\pm 4932}$ & $2^{-113} \approx 9.6\cdot 10^{-35}$
\end{tabular}
\end{table}

\subsection{Experiments with $(u_f, u, u_r) = $ (double, double, quad)}
Table \ref{tab:sdq} shows the experiments for the setting $(u_f, u, u_r)=$ (double, double, quad), with both ${\epsilon}=2^{-53}$ and ${\epsilon}=2^{-37}$.
First, we note that where SPAI-GMRES-IR converges, BSPAI-GMRES-IR also converges, as predicted by our theoretical results, although of course the adaptive precision storage can result in a different total number of GMRES iterations across the refinement steps.  

The performance for the matrix \texttt{steam1} using ${\varepsilon}=0.1$ with  ${\epsilon}=2^{-53}$ and ${\epsilon}=2^{-37}$, BSPAI-GMRES-IR takes $21$ total GMRES iterations to converge to double precision accuracy while SPAI-GMRES-IR takes total $14$ GMRES iterations. The storage (and computation) savings of the adaptive precision approach can be significant for this case; using $\epsilon = 2^{-53}$ and $\epsilon=2^{-37}$, requires only 74.9\% and 42.6\% of the storage/computation cost as the uniform precision approach, respectively. This matrix perhaps represents a best-case scenario. For \texttt{steam3}, we also see reasonable reductions in storage cost for the two choices of $\epsilon$; note that although the choice $\epsilon=2^{-37}$ results in significant storage savings, the number of GMRES iterations required increases significantly. For some matrices, such as \texttt{bfwa782} and \texttt{sherman4}, there appears to be no benefit to the adaptive precision approach. 

\begin{table}[]
\caption{Comparison of BSPAI-GMRES-IR for different ${\epsilon}$ values with SPAI-GMRES-IR for test matrices using $(u_{f},u,u_{r})=$ (double, double, quad). Note that the $\varepsilon$ parameter indicated for SPAI is also used for BSPAI.} 
\begin{center}
     {\scriptsize{
\begin{tabular}{|c|c|c|c|c|c|}
\hline
& Preconditioner & $\kappa_\infty(\tilde{A})$ & Precond. $nnz$ & str.& \begin{tabular}[c]{@{}c@{}}GMRES its/step\end{tabular} \\ \hline

\parbox[t]{9mm}{\multirow{3}{*}{\rotatebox[origin=c]{0}{\texttt{steam1}}}} & 
BSPAI(${\epsilon}=2^{-53}$)& {\multirow{3}{*}{\rotatebox[origin=c]{0}1.5}} &  1105(556, 537, 12, 0) &  74.9\% & 21(7, 7, 7) \\ 
& BSPAI(${\epsilon}=2^{-37}$) & 
& 1105(242, 284, 347, 232) &  42.6\% & 21(7, 7, 7) \\ 
& SPAI (${\varepsilon}=0.1$) &  & 1105(1105, 0, 0,0) &  &14(7, 7) \\ \hline

\parbox[t]{9mm}{\multirow{3}{*}{\rotatebox[origin=c]{0}{\texttt{pores\_3}}}} &  
BSPAI (${\epsilon}=2^{-53}$) & 
  & 3560(3404, 156, 0, 0) & 97.8\% &323(111, 103, 109) \\
& BSPAI (${\epsilon}=2^{-37}$) & $6.6\cdot10^{3}$ & 3560(2409, 924, 218, 9) & 82.1\%  & 333(111, 112, 110) \\ 
& SPAI (${\varepsilon}=0.5$) &  & 3560(3560, 0, 0,0) &  & 223(111, 112) \\ \hline

\parbox[t]{9mm}{\multirow{3}{*}{\rotatebox[origin=c]{0}{\texttt{steam3}}}} 
&BSPAI (${\epsilon}=2^{-53}$) &  & 347(248, 83, 14, 2) &  84.4\% & 15(5, 5, 5) \\ 
&BSPAI (${\epsilon}=2^{-37}$) & $1.9$ & 347(139, 85, 92, 31) & 59.0\% & 80(80) \\ 
&SPAI ($\varepsilon=0.1$) &  & 347(347, 0, 0, 0) &  & 11(5, 6) \\ \hline

\parbox[t]{9mm}{\multirow{3}{*}{\rotatebox[origin=c]{0}{\texttt{saylr1}}}} & 
BSPAI (${\epsilon}=2^{-53}$) &  & 1932(1890, 42, 0, 0) & 98.9\% & 197(64, 66, 67) \\ 
& BSPAI (${\epsilon}=2^{-37}$) &  $1.9\cdot 10^4$ & 1932(977, 830, 125, 0) & 73.7\% & 238(238) \\ 
& SPAI ($\varepsilon=0.4$) &  & 1932(1932, 0, 0, 0) &  & 195(64, 66, 65) \\ \hline
\parbox[t]{9mm}{\multirow{3}{*}{\rotatebox[origin=c]{0}{\texttt{cage5}}}} & 
BSPAI (${\epsilon}=2^{-53}$) &  &  511(511, 0, 0, 0) & 100\% & 14(7, 7) \\
& BSPAI (${\epsilon}=2^{-37}$) &  $1.8$ & 511(504, 7, 0, 0) & 99.3\% & 14(7, 7) \\
& SPAI (${\varepsilon}=0.1$) & & 511(511, 0, 0, 0) &  & 14(7, 7) \\ \hline

\parbox[t]{9mm}{\multirow{3}{*}{\rotatebox[origin=c]{0}{\texttt{gre\_115}}}} &   
BSPAI (${\epsilon}=2^{-53}$) & & 5095(5094, 1, 0, 0) &  100\% & 18(9, 9) \\ 
& BSPAI (${\epsilon}=2^{-37}$) & $3.2$ & 5095(4145, 949, 1, 0) & 90.7\% & 27(9, 9, 9) \\ 
& SPAI (${\varepsilon}=0.1$) & & 5095(5095, 0, 0, 0) & & 18(9, 9) \\ \hline

\parbox[t]{9mm}{\multirow{3}{*}{\rotatebox[origin=c]{0}{\hspace{-1.5mm}\texttt{sherman4}}}} &   
BSPAI (${\epsilon}=2^{(-53)}$) &  &  1385(1385, 0, 0, 0) & 100\% & 187(92, 95) \\ 
& BSPAI (${\epsilon}=2^{(-37)}$) & $1.6\cdot 10^{3}$ & 1385(1382, 3, 0, 0) & 99.9\%  & 276(92, 92, 92) \\ 
& SPAI (${\varepsilon}=0.5$) &  & 1385(1385, 0, 0, 0) &  & 187(92, 95) \\ \hline

\parbox[t]{9mm}{\multirow{3}{*}{\rotatebox[origin=c]{0}{\texttt{hor\_131}}}} & 
BSPAI (${\epsilon}=2^{-53}$) &  & 8478(8471, 7, 0, 0) & 100\%  & \hspace{-1mm}436(108, 111, 109, 108)\hspace{-1mm} \\ 
& BSPAI (${\epsilon}=2^{-37}$) & $1.3\cdot 10^5$ & 8478(4736, 3707, 35, 0) & 77.8\% & \hspace{-1mm}440(108, 110, 111, 111)\hspace{-1mm} \\ 
& SPAI (${\varepsilon}=0.5$) &  & 8478(8478, 0, 0,0) &  & 327(108, 111, 108) \\ \hline

\parbox[t]{9mm}{\multirow{3}{*}{\rotatebox[origin=c]{0}{\texttt{bfwa782}}}} &
BSPAI (${\epsilon}=2^{-53}$) &  & 6261(6261, 0, 0, 0) & 100\% & 234(113, 121) \\ 
& BSPAI (${\epsilon}=2^{-37}$) &  $1.1\cdot 10^3$ & 6261(6040, 221, 0, 0) & 98.2\% & 355(113, 121, 121) \\ 
& SPAI (${\varepsilon}=0.5$ &  & 6261(6261, 0, 0, 0) & & 234(113, 121) \\ \hline

\end{tabular}
}}
\end{center}
\label{tab:sdq}
\end{table}

\subsection{Experiments with $(u_f, u, u_r) = $ (single, single, double)}
Table \ref{tab:ssd} shows the experiments for the setting  $(u_f, u, u_r)=$ (single, single, double) with $u_1$ = single, $u_2=$ half, and $u_3=$ 1, and with $\epsilon=2^{-24}$ and ${\epsilon}=2^{-18}$. In all tests, both SPAI-GMRES-IR and BSPAI-GMRES-IR converge to single precision accuracy. 

For the value $\epsilon=2^{-24}$, BSPAI-GMRES-IR  takes about the same number of iterations as SPAI-GMRES-IR and requires an average of 98.6\% of the storage cost as the uniform precision approach. The performance of BSPAI-GMRES-IR for the matrix \texttt{cage5} using $\epsilon=2^{-18}$ requires 76.5\% of the storage cost as uniform precision and converges in the same number of iteration as SPAI-GMRES-IR. For matrices \texttt{gre\_115}, \texttt{sherman4}, and \texttt{bfwa782}, although using $\epsilon=2^{-18}$ results in storage costs of 66.7\%, 70.7\% and 73.8\% of the uniform precision approach, respectively, a greater number of iterations are required than in the uniform precision SPAI-GMRES-IR. 

\begin{table}[h!]
\caption{Comparison of BSPAI-GMRES-IR for different ${\epsilon}$ values with SPAI-GMRES-IR for test matrices using $(u_{f},u,u_{r})=$ (single, single, double). Note that the $\varepsilon$ parameter indicated for SPAI is also used for BSPAI.}
\begin{center}
     {\scriptsize{
\begin{tabular}{|c|c|c|c|c|c|}
\hline
& Preconditioner & $\kappa_\infty(\tilde{A})$ & Precond. $nnz$ & str. & \begin{tabular}[c]{@{}c@{}}GMRES its/step\end{tabular} \\ \hline

\parbox[t]{9mm}{\multirow{3}{*}{\rotatebox[origin=c]{0}{\texttt{cage5}}}} & 
BSPAI (${\epsilon}=2^{-24}$) & {\multirow{3}{*}{\rotatebox[origin=c]{0}{$1.8$}}} &  511(504, 7, 0) & 99.3\% & 8(4, 4) \\
& BSPAI (${\epsilon}=2^{-18}$) & &  511(271, 239, 1) & 76.5\% & 8(4, 4) \\
& SPAI (${\varepsilon}=0.1$) &  & 511(511, 0, 0) &  & 8(4, 4) \\ \hline

\parbox[t]{9mm}{\multirow{3}{*}{\rotatebox[origin=c]{0}{\texttt{gre\_115}}}} &   
BSPAI (${\epsilon}=2^{-24}$) & {\multirow{3}{*}{\rotatebox[origin=c]{0}{$6.6\cdot 10^1$}}} & 1094(1052, 40, 2) & 98\% & 32(16, 16) \\ 
& BSPAI (${\epsilon}=2^{-18}$) & & 1094(367, 724, 3) & 66.7\% & 50(16, 17, 17) \\ 
& SPAI (${\varepsilon}=0.4$) &  & 1094(1094, 0, 0) &  & 32(16, 16) \\ \hline

\parbox[t]{9mm}{\multirow{3}{*}{\rotatebox[origin=c]{0}{\hspace{-1.5mm}\texttt{sherman4}}}} &   
BSPAI (${\epsilon}=2^{(-24)}$) & {\multirow{3}{*}{\rotatebox[origin=c]{0}{$5.0\cdot 10^2$}}}&  8499(8289, 210, 0) & 98.8\% & 73(35, 38) \\ 
& BSPAI (${\epsilon}=2^{(-18)}$) &  &  8499(3527, 4957, 15) & 70.7\% & 103(35, 35, 33) \\ 
& SPAI (${\varepsilon}=0.3$) & & 8499(8499, 0, 0) &  & 74(35, 39) \\ \hline

\parbox[t]{9mm}{\multirow{3}{*}{\rotatebox[origin=c]{0}{\texttt{bfwa782}}}} &
BSPAI (${\epsilon}=2^{-24}$) & {\multirow{3}{*}{\rotatebox[origin=c]{0}{$6.6\cdot 10^2$}}} & 8053(7774, 279, 0) & 98.3\% & 127(57, 70) \\ 
& BSPAI (${\epsilon}=2^{-18}$) & & 8053(3832, 4217, 4) & 73.8\% & 191(57, 68, 66) \\ 
& SPAI (${\varepsilon}=0.4$) &  & 8053(8053, 0, 0) &  & 127(57, 70) \\ \hline

\parbox[t]{9mm}{\multirow{2}{*}{\rotatebox[origin=c]{0}{\texttt{pores\_3}}}} & 
BSPAI (${\epsilon}=2^{-24}$) & {\multirow{2}{*}{\rotatebox[origin=c]{0}{$6.6 \cdot 10^3$}}} &  3560(2403, 876, 281) & 81.8\% & 353(93  93  90  77) \\ 
& SPAI (${\varepsilon}=0.5$) & & 3560(3560, 0, 0) &  & 188(93, 95) \\ \hline

\parbox[t]{9mm}{\multirow{2}{*}{\rotatebox[origin=c]{0}{\texttt{steam1}}}} & 
BSPAI (${\epsilon}=2^{-24}$) & {\multirow{2}{*}{\rotatebox[origin=c]{0}{$1.5$}}} &  1302(259,  295,  748) & 45.6\% & 244(4, 240) \\ 
& SPAI (${\varepsilon}=0.1$) &  & 1302(1302, 0, 0) &  & 8(4, 4) \\ \hline

\parbox[t]{9mm}{\multirow{2}{*}{\rotatebox[origin=c]{0}{\texttt{steam3}}}} & 
BSPAI (${\epsilon}=2^{-24}$) & {\multirow{2}{*}{\rotatebox[origin=c]{0}{$1.9$}}} &  429(135, 85, 209) & 53.6\% & 83(3, 80) \\ 

& SPAI (${\varepsilon}=0.1$) &  & 429(429, 0, 0) &  & 6(3, 3) \\ \hline

\parbox[t]{9mm}{\multirow{2}{*}{\rotatebox[origin=c]{0}{\texttt{saylr1}}}} & 
BSPAI (${\epsilon}=2^{-24}$) & {\multirow{2}{*}{\rotatebox[origin=c]{0}{$1.9 \cdot 10^4$}}} &  1932(977, 773, 182) & 72.9\% & 238(238) \\ 
& SPAI (${\varepsilon}=0.4$) &  & 1932(1932, 0, 0) &  & 154(40, 57, 57) \\ \hline
\end{tabular}
}}
\end{center}
\label{tab:ssd}
\end{table}

\section{Conclusions and Future Work}
\label{sec:conclusion}
In this work we use an adaptive precision sparse approximate inverse preconditioner within mixed precision GMRES-based iterative refinement. Using the approach of Graillat et al. \cite{graillat2022adaptive}, after computing a sparse approximate inverse in low precision, we place elements of the sparse approximate inverse preconditioner into buckets for a given set of precisions based on their magnitude. We then apply the preconditioner to a vector in mixed precision within five precision GMRES-IR; we call this algorithm variant BSPAI-GMRES-IR. 

We then analyze the behavior of the backward and forward errors of mixed precision left-preconditioned GMRES method,  which uses the bucketed sparse approximate inverse as a left preconditioner. Our analysis shows that if we choose $u_1\approx \varepsilon \approx u_p$, then the normwise backward and forward errors will be close to those we get in the case that we use uniform precision. This indicates that BSPAI-GMRES-IR will converge under the same conditions as SPAI-GMRES-IR. 

We performing a set of numerical experiments which shows that the adaptive sparse-matrix vector product approach can reduce the cost of storing and applying the sparse approximate inverse preconditioner, although a significant reduction in cost often comes at the expense of increasing the number of GMRES iterations required for convergence.  
We note that it is possible to extend this approach to other preconditioners for Krylov subspace methods.

We again stress that a fruitful potential area of future work is to extend the adaptive sparse-matrix vector product approach to decouple the storage and computation precisions. This would make this approach beneficial for existing cases where we would ideally like to store a matrix in lower precision and apply it to a vector in a higher precision, which is often the case within SPAI-GMRES-IR. 

Other potential future work involves the development and analysis of other adaptive-precision matrix computations, such as triangular solves.

\bibliographystyle{siamplain}
\bibliography{BSPAI}
\end{document}